\documentclass[12pt]{amsart}
\usepackage{epsfig}

\makeatletter
 
\def\diagram{\m@th\leftwidth=\z@ \rightwidth=\z@ \topheight=\z@
\botheight=\z@ \setbox\@picbox\hbox\bgroup}
 
\def\enddiagram{\egroup\wd\@picbox\rightwidth\unitlength
\ht\@picbox\topheight\unitlength \dp\@picbox\botheight\unitlength
\hskip\leftwidth\unitlength\box\@picbox}
 
\def\bfig{\begin{diagram}}
\def\efig{\end{diagram}}
\newcount\wideness \newcount\leftwidth \newcount\rightwidth
\newcount\highness \newcount\topheight \newcount\botheight
 
\def\ratchet#1#2{\ifnum#1<#2 \global #1=#2 \fi}
 
\def\putbox(#1,#2)#3{%
\horsize{\wideness}{#3} \divide\wideness by 2
{\advance\wideness by #1 \ratchet{\rightwidth}{\wideness}}
{\advance\wideness by -#1 \ratchet{\leftwidth}{\wideness}}
\vertsize{\highness}{#3} \divide\highness by 2
{\advance\highness by #2 \ratchet{\topheight}{\highness}}
{\advance\highness by -#2 \ratchet{\botheight}{\highness}}
\put(#1,#2){\makebox(0,0){$#3$}}}
 
\def\putlbox(#1,#2)#3{%
\horsize{\wideness}{#3}
{\advance\wideness by #1 \ratchet{\rightwidth}{\wideness}}
{\ratchet{\leftwidth}{-#1}}
\vertsize{\highness}{#3} \divide\highness by 2
{\advance\highness by #2 \ratchet{\topheight}{\highness}}
{\advance\highness by -#2 \ratchet{\botheight}{\highness}}
\put(#1,#2){\makebox(0,0)[l]{$#3$}}}
 
\def\putrbox(#1,#2)#3{%
\horsize{\wideness}{#3}
{\ratchet{\rightwidth}{#1}}
{\advance\wideness by -#1 \ratchet{\leftwidth}{\wideness}}
\vertsize{\highness}{#3} \divide\highness by 2
{\advance\highness by #2 \ratchet{\topheight}{\highness}}
{\advance\highness by -#2 \ratchet{\botheight}{\highness}}
\put(#1,#2){\makebox(0,0)[r]{$#3$}}}

\def\adjust[#1]{} 
 
\newcount \coefa
\newcount \coefb
\newcount \coefc
\newcount\tempcounta
\newcount\tempcountb
\newcount\tempcountc
\newcount\tempcountd
\newcount\xext
\newcount\yext
\newcount\xoff
\newcount\yoff
\newcount\gap%
\newcount\arrowtypea
\newcount\arrowtypeb
\newcount\arrowtypec
\newcount\arrowtyped
\newcount\arrowtypee
\newcount\height
\newcount\width
\newcount\xpos
\newcount\ypos
\newcount\run
\newcount\rise
\newcount\arrowlength
\newcount\halflength
\newcount\arrowtype
\newdimen\tempdimen
\newdimen\xlen
\newdimen\ylen
\newsavebox{\tempboxa}%
\newsavebox{\tempboxb}%
\newsavebox{\tempboxc}%
 
\newdimen\w@dth
 
\def\setw@dth#1#2{\setbox\z@\hbox{\m@th$#1$}\w@dth=\wd\z@
\setbox\@ne\hbox{\m@th$#2$}\ifnum\w@dth<\wd\@ne \w@dth=\wd\@ne \fi
\advance\w@dth by 1.2em}
 
\def\t@^#1_#2{\allowbreak\def\n@one{#1}\def\n@two{#2}\mathrel
{\setw@dth{#1}{#2}
\mathop{\hbox to \w@dth{\rightarrowfill}}\limits
\ifx\n@one\empty\else ^{\box\z@}\fi
\ifx\n@two\empty\else _{\box\@ne}\fi}}
\def\t@@^#1{\@ifnextchar_{\t@^{#1}}{\t@^{#1}_{}}}
\def\to{\@ifnextchar^{\t@@}{\t@@^{}}}
 
\def\t@left^#1_#2{\def\n@one{#1}\def\n@two{#2}\mathrel{\setw@dth{#1}{#2}
\mathop{\hbox to \w@dth{\leftarrowfill}}\limits
\ifx\n@one\empty\else ^{\box\z@}\fi
\ifx\n@two\empty\else _{\box\@ne}\fi}}
\def\t@@left^#1{\@ifnextchar_{\t@left^{#1}}{\t@left^{#1}_{}}}
\def\toleft{\@ifnextchar^{\t@@left}{\t@@left^{}}}
 
\def\two@^#1_#2{\allowbreak
\def\n@one{#1}\def\n@two{#2}\mathrel{\setw@dth{#1}{#2}
\mathop{\vcenter{\lineskip\z@\baselineskip\z@
                 \hbox to \w@dth{\rightarrowfill}%
                 \hbox to \w@dth{\rightarrowfill}}%
       }\limits
\ifx\n@one\empty\else ^{\box\z@}\fi
\ifx\n@two\empty\else _{\box\@ne}\fi}}
\def\tw@@^#1{\@ifnextchar _{\two@^{#1}}{\two@^{#1}_{}}}
\def\two{\@ifnextchar ^{\tw@@}{\tw@@^{}}}
 
\def\tofr@^#1_#2{\def\n@one{#1}\def\n@two{#2}\mathrel{\setw@dth{#1}{#2}
\mathop{\vcenter{\hbox to \w@dth{\rightarrowfill}\kern-1.7ex
                 \hbox to \w@dth{\leftarrowfill}}%
       }\limits
\ifx\n@one\empty\else ^{\box\z@}\fi
\ifx\n@two\empty\else _{\box\@ne}\fi}}
\def\t@fr@^#1{\@ifnextchar_ {\tofr@^{#1}}{\tofr@^{#1}_{}}}
\def\tofro{\@ifnextchar^ {\t@fr@}{\t@fr@^{}}}

\def\mon{\mathop{\m@th\hbox to
      14.6\P@{\lasyb\char'51\hskip-2.1\P@$\arrext$\hss
$\mathord\rightarrow$}}\limits} 
\def\leftmono{\mathrel{\m@th\hbox to
14.6\P@{$\mathord\leftarrow$\hss$\arrext$\hskip-2.1\P@\lasyb\char'50%
}}\limits} 
\mathchardef\arrext="0200       

\def\settypes(#1,#2,#3){\arrowtypea#1 \arrowtypeb#2 \arrowtypec#3}
\def\settoheight#1#2{\setbox\@tempboxa\hbox{#2}#1\ht\@tempboxa\relax}%
\def\settodepth#1#2{\setbox\@tempboxa\hbox{#2}#1\dp\@tempboxa\relax}%
\def\settokens`#1`#2`#3`#4`{%
     \def\tokena{#1}\def\tokenb{#2}\def\tokenc{#3}\def\tokend{#4}}
\def\setsqparms[#1`#2`#3`#4;#5`#6]{%
\arrowtypea #1
\arrowtypeb #2
\arrowtypec #3
\arrowtyped #4
\width #5
\height #6
}
\def\setpos(#1,#2){\xpos=#1 \ypos#2}

\def\settriparms[#1`#2`#3;#4]{\settripairparms[#1`#2`#3`1`1;#4]}%
 
\def\settripairparms[#1`#2`#3`#4`#5;#6]{%
\arrowtypea #1
\arrowtypeb #2
\arrowtypec #3
\arrowtyped #4
\arrowtypee #5
\width #6
\height #6
}
 
\def\resetparms{\settripairparms[1`1`1`1`1;500]\width 500}
 
\resetparms
 
\def\mvector(#1,#2)#3{
\put(0,0){\vector(#1,#2){#3}}%
\put(0,0){\vector(#1,#2){26}}%
}
\def\evector(#1,#2)#3{{
\arrowlength #3
\put(0,0){\vector(#1,#2){\arrowlength}}%
\advance \arrowlength by-30
\put(0,0){\vector(#1,#2){\arrowlength}}%
}}
 
\def\horsize#1#2{%
\settowidth{\tempdimen}{$#2$}%
#1=\tempdimen
\divide #1 by\unitlength
}
 
\def\vertsize#1#2{%
\settoheight{\tempdimen}{$#2$}%
#1=\tempdimen
\settodepth{\tempdimen}{$#2$}%
\advance #1 by\tempdimen
\divide #1 by\unitlength
}
 
\def\putvector(#1,#2)(#3,#4)#5#6{{%
\ifnum3<\arrowtype
\putdashvector(#1,#2)(#3,#4)#5\arrowtype
\else
\ifnum\arrowtype<-3
\putdashvector(#1,#2)(#3,#4)#5\arrowtype
\else
\xpos=#1
\ypos=#2
\run=#3
\rise=#4
\arrowlength=#5
\ifnum \arrowtype<0
    \ifnum \run=0
        \advance \ypos by-\arrowlength
    \else
        \tempcounta \arrowlength
        \multiply \tempcounta by\rise
        \divide \tempcounta by\run
        \ifnum\run>0
            \advance \xpos by\arrowlength
            \advance \ypos by\tempcounta
        \else
            \advance \xpos by-\arrowlength
            \advance \ypos by-\tempcounta
        \fi
    \fi
    \multiply \arrowtype by-1
    \multiply \rise by-1
    \multiply \run by-1
\fi
\ifcase \arrowtype
\or \put(\xpos,\ypos){\vector(\run,\rise){\arrowlength}}%
\or \put(\xpos,\ypos){\mvector(\run,\rise)\arrowlength}%
\or \put(\xpos,\ypos){\evector(\run,\rise){\arrowlength}}%
\fi\fi\fi
}}
 
\def\putsplitvector(#1,#2)#3#4{
\xpos #1
\ypos #2
\arrowtype #4
\halflength #3
\arrowlength #3
\gap 140
\advance \halflength by-\gap
\divide \halflength by2
\ifnum\arrowtype>0
   \ifcase \arrowtype
   \or \put(\xpos,\ypos){\line(0,-1){\halflength}}%
       \advance\ypos by-\halflength
       \advance\ypos by-\gap
       \put(\xpos,\ypos){\vector(0,-1){\halflength}}%
   \or \put(\xpos,\ypos){\line(0,-1)\halflength}%
       \put(\xpos,\ypos){\vector(0,-1)3}%
       \advance\ypos by-\halflength
       \advance\ypos by-\gap
       \put(\xpos,\ypos){\vector(0,-1){\halflength}}%
   \or \put(\xpos,\ypos){\line(0,-1)\halflength}%
       \advance\ypos by-\halflength
       \advance\ypos by-\gap
       \put(\xpos,\ypos){\evector(0,-1){\halflength}}%
   \fi
\else \arrowtype=-\arrowtype
   \ifcase\arrowtype
   \or \advance \ypos by-\arrowlength
       \put(\xpos,\ypos){\line(0,1){\halflength}}%
       \advance\ypos by\halflength
       \advance\ypos by\gap
       \put(\xpos,\ypos){\vector(0,1){\halflength}}%
   \or \advance \ypos by-\arrowlength
       \put(\xpos,\ypos){\line(0,1)\halflength}%
       \put(\xpos,\ypos){\vector(0,1)3}%
       \advance\ypos by\halflength
       \advance\ypos by\gap
       \put(\xpos,\ypos){\vector(0,1){\halflength}}%
   \or \advance \ypos by-\arrowlength
       \put(\xpos,\ypos){\line(0,1)\halflength}%
       \advance\ypos by\halflength
       \advance\ypos by\gap
       \put(\xpos,\ypos){\evector(0,1){\halflength}}%
   \fi
\fi
}
 
\def\putmorphism(#1)(#2,#3)[#4`#5`#6]#7#8#9{{%
\run #2
\rise #3
\ifnum\rise=0
  \puthmorphism(#1)[#4`#5`#6]{#7}{#8}#9%
\else\ifnum\run=0
  \putvmorphism(#1)[#4`#5`#6]{#7}{#8}#9%
\else
\setpos(#1)%
\arrowlength #7
\arrowtype #8
\ifnum\run=0
\else\ifnum\rise=0
\else
\ifnum\run>0
    \coefa=1
\else
   \coefa=-1
\fi
\ifnum\arrowtype>0
   \coefb=0
   \coefc=-1
\else
   \coefb=\coefa
   \coefc=1
   \arrowtype=-\arrowtype
\fi
\width=2
\multiply \width by\run
\divide \width by\rise
\ifnum \width<0  \width=-\width\fi
\advance\width by60
\if l#9 \width=-\width\fi
\putbox(\xpos,\ypos){#4}
{\multiply \coefa by\arrowlength
\advance\xpos by\coefa
\multiply \coefa by\rise
\divide \coefa by\run
\advance \ypos by\coefa
\putbox(\xpos,\ypos){#5} }%
{\multiply \coefa by\arrowlength
\divide \coefa by2
\advance \xpos by\coefa
\advance \xpos by\width
\multiply \coefa by\rise
\divide \coefa by\run
\advance \ypos by\coefa
\if l#9%
   \putrbox(\xpos,\ypos){#6}%
\else\if r#9%
   \putlbox(\xpos,\ypos){#6}%
\fi\fi }%
{\multiply \rise by-\coefc
\multiply \run by-\coefc
\multiply \coefb by\arrowlength
\advance \xpos by\coefb
\multiply \coefb by\rise
\divide \coefb by\run
\advance \ypos by\coefb
\multiply \coefc by70
\advance \ypos by\coefc
\multiply \coefc by\run
\divide \coefc by\rise
\advance \xpos by\coefc
\multiply \coefa by140
\multiply \coefa by\run
\divide \coefa by\rise
\advance \arrowlength by\coefa
\ifcase\arrowtype
\or \put(\xpos,\ypos){\vector(\run,\rise){\arrowlength}}%
\or \put(\xpos,\ypos){\mvector(\run,\rise){\arrowlength}}%
\or \put(\xpos,\ypos){\evector(\run,\rise){\arrowlength}}%
\fi}\fi\fi\fi\fi}}

\newcount\numbdashes \newcount\lengthdash \newcount\increment
 
\def\howmanydashes{
\numbdashes=\arrowlength \lengthdash=40
\divide\numbdashes by \lengthdash
\lengthdash=\arrowlength
\divide\lengthdash by \numbdashes
\increment=\lengthdash
\multiply\lengthdash by 3
\divide\lengthdash by 5
}
 
\def\putdashvector(#1)(#2,#3)#4#5{%
\ifnum#3=0 \putdashhvector(#1){#4}#5
\else
\ifnum#2=0
\putdashvvector(#1){#4}#5\fi\fi}
 
\def\putdashhvector(#1,#2)#3#4{{%
\arrowlength=#3 \howmanydashes
\multiput(#1,#2)(\increment,0){\numbdashes}%
{\vrule height .4pt width \lengthdash\unitlength}
\arrowtype=#4 \xpos=#1
\ifnum\arrowtype<0 \advance\arrowtype by 7 \fi
\ifcase\arrowtype
\or \advance\xpos by 10
    \put(\xpos,#2){\vector(-1,0){\lengthdash}}
    \advance\xpos by 40
    \put(\xpos,#2){\vector(-1,0){\lengthdash}}
\or \advance \xpos by 10
    \put(\xpos,#2){\vector(-1,0){\lengthdash}}
    \advance\xpos by  \arrowlength
    \advance\xpos by  -50
    \put(\xpos,#2){\vector(-1,0){\lengthdash}}
\or \advance\xpos by 10
    \put(\xpos,#2){\vector(-1,0){\lengthdash}}
\or \advance\xpos by \arrowlength
    \advance\xpos by -\lengthdash
    \put(\xpos,#2){\vector(1,0){\lengthdash}}
\or {\advance\xpos by 10
    \put(\xpos,#2){\vector(1,0){\lengthdash}}}
    \advance\xpos by \arrowlength
    \advance\xpos by -\lengthdash
    \put(\xpos,#2){\vector(1,0){\lengthdash}}
\or \advance\xpos by \arrowlength
    \advance\xpos by -\lengthdash
    \put(\xpos,#2){\vector(1,0){\lengthdash}}
    \advance\xpos by -40
    \put(\xpos,#2){\vector(1,0){\lengthdash}}
   \fi
}}
 
\def\putdashvvector(#1,#2)#3#4{{%
\arrowlength=#3 \howmanydashes
\ypos=#2 \advance\ypos by -\arrowlength
\multiput(#1,#2)(0,\increment){\numbdashes}%
    {\vrule width .4pt height \lengthdash\unitlength}
\arrowtype=#4 \ypos=#2
\ifnum\arrowtype<0 \advance\arrowtype by 7 \fi
\ifcase\arrowtype
\or \advance\ypos by \arrowlength \advance\ypos by -40
    \put(#1,\ypos){\vector(0,1){\lengthdash}}
    \advance\ypos by -40
    \put(#1,\ypos){\vector(0,1){\lengthdash}}
\or \advance\ypos by 10
    \put(#1,\ypos){\vector(0,1){\lengthdash}}
    \advance\ypos by \arrowlength \advance\ypos by -40
    \put(#1,\ypos){\vector(0,1){\lengthdash}}
\or \advance\ypos by \arrowlength \advance\ypos by -40
    \put(#1,\ypos){\vector(0,1){\lengthdash}}
\or \advance\ypos by 10
    \put(#1,\ypos){\vector(0,-1){\lengthdash}}
\or \advance\ypos by 10
    \put(#1,\ypos){\vector(0,-1){\lengthdash}}
    \advance\ypos by \arrowlength \advance\ypos by -40
    \put(#1,\ypos){\vector(0,-1){\lengthdash}}
\or \advance\ypos by 10
    \put(#1,\ypos){\vector(0,-1){\lengthdash}}
    \advance\ypos by 40
    \put(#1,\ypos){\vector(0,-1){\lengthdash}}
\fi
}}
 
\def\puthmorphism(#1,#2)[#3`#4`#5]#6#7#8{{%
\xpos #1
\ypos #2
\width #6
\arrowlength #6
\arrowtype=#7
\putbox(\xpos,\ypos){#3\vphantom{#4}}%
{\advance \xpos by\arrowlength
\putbox(\xpos,\ypos){\vphantom{#3}#4}}%
\horsize{\tempcounta}{#3}%
\horsize{\tempcountb}{#4}%
\divide \tempcounta by2
\divide \tempcountb by2
\advance \tempcounta by30
\advance \tempcountb by30
\advance \xpos by\tempcounta
\advance \arrowlength by-\tempcounta
\advance \arrowlength by-\tempcountb
\putvector(\xpos,\ypos)(1,0)\arrowlength\arrowtype
\divide \arrowlength by2
\advance \xpos by\arrowlength
\vertsize{\tempcounta}{#5}%
\divide\tempcounta by2
\advance \tempcounta by20
\if a#8 %
   \advance \ypos by\tempcounta
   \putbox(\xpos,\ypos){#5}%
\else
   \advance \ypos by-\tempcounta
   \putbox(\xpos,\ypos){#5}%
\fi}}
 
\def\putvmorphism(#1,#2)[#3`#4`#5]#6#7#8{{%
\xpos #1
\ypos #2
\arrowlength #6
\arrowtype #7
\settowidth{\xlen}{$#5$}%
\putbox(\xpos,\ypos){#3}%
{\advance \ypos by-\arrowlength
\putbox(\xpos,\ypos){#4}}%
{\advance\arrowlength by-140
\advance \ypos by-70
\ifdim\xlen>0pt
   \if m#8%
      \putsplitvector(\xpos,\ypos)\arrowlength\arrowtype
   \else
   \putvector(\xpos,\ypos)(0,-1)\arrowlength\arrowtype
   \fi
\else
   \putvector(\xpos,\ypos)(0,-1)\arrowlength\arrowtype
\fi}%
\ifdim\xlen>0pt
   \divide \arrowlength by2
   \advance\ypos by-\arrowlength
   \if l#8%
      \advance \xpos by-40
      \putrbox(\xpos,\ypos){#5}%
   \else\if r#8%
      \advance \xpos by40
      \putlbox(\xpos,\ypos){#5}%
   \else
      \putbox(\xpos,\ypos){#5}%
   \fi\fi
\fi
}}
 
\def\putsquarep<#1>(#2)[#3;#4`#5`#6`#7]{{%
\setsqparms[#1]%
\setpos(#2)%
\settokens`#3`%
\puthmorphism(\xpos,\ypos)[\tokenc`\tokend`{#7}]{\width}{\arrowtyped}b%
\advance\ypos by \height
\puthmorphism(\xpos,\ypos)[\tokena`\tokenb`{#4}]{\width}{\arrowtypea}a%
\putvmorphism(\xpos,\ypos)[``{#5}]{\height}{\arrowtypeb}l%
\advance\xpos by \width
\putvmorphism(\xpos,\ypos)[``{#6}]{\height}{\arrowtypec}r%
}}
 
\def\putsquare{\@ifnextchar <{\putsquarep}{\putsquarep%
   <\arrowtypea`\arrowtypeb`\arrowtypec`\arrowtyped;\width`\height>}}
\def\square{\@ifnextchar< {\squarep}{\squarep
   <\arrowtypea`\arrowtypeb`\arrowtypec`\arrowtyped;\width`\height>}}
\def\squarep<#1>[#2`#3`#4`#5;#6`#7`#8`#9]{{
\setsqparms[#1]
\diagram
\putsquarep<\arrowtypea`\arrowtypeb`\arrowtypec`
\arrowtyped;\width`\height>
(0,0)[#2`#3`#4`{#5};#6`#7`#8`{#9}]
\enddiagram
}}                                                 
\def\putptrianglep<#1>(#2,#3)[#4`#5`#6;#7`#8`#9]{{%
\settriparms[#1]%
\xpos=#2 \ypos=#3
\advance\ypos by \height
\puthmorphism(\xpos,\ypos)[#4`#5`{#7}]{\height}{\arrowtypea}a%
\putvmorphism(\xpos,\ypos)[`#6`{#8}]{\height}{\arrowtypeb}l%
\advance\xpos by\height
\putmorphism(\xpos,\ypos)(-1,-1)[``{#9}]{\height}{\arrowtypec}r%
}}
 
\def\putptriangle{\@ifnextchar <{\putptrianglep}{\putptrianglep
   <\arrowtypea`\arrowtypeb`\arrowtypec;\height>}}
\def\ptriangle{\@ifnextchar <{\ptrianglep}{\ptrianglep
   <\arrowtypea`\arrowtypeb`\arrowtypec;\height>}}
\def\ptrianglep<#1>[#2`#3`#4;#5`#6`#7]{{
\settriparms[#1]
\diagram
\putptrianglep<\arrowtypea`\arrowtypeb`
\arrowtypec;\height>
(0,0)[#2`#3`#4;#5`#6`{#7}]
\enddiagram
}}                                            
 
\def\putqtrianglep<#1>(#2,#3)[#4`#5`#6;#7`#8`#9]{{%
\settriparms[#1]%
\xpos=#2 \ypos=#3
\advance\ypos by\height
\puthmorphism(\xpos,\ypos)[#4`#5`{#7}]{\height}{\arrowtypea}a%
\putmorphism(\xpos,\ypos)(1,-1)[``{#8}]{\height}{\arrowtypeb}l%
\advance\xpos by\height
\putvmorphism(\xpos,\ypos)[`#6`{#9}]{\height}{\arrowtypec}r%
}}
 
\def\putqtriangle{\@ifnextchar <{\putqtrianglep}{\putqtrianglep
   <\arrowtypea`\arrowtypeb`\arrowtypec;\height>}}
\def\qtriangle{\@ifnextchar <{\qtrianglep}{\qtrianglep
   <\arrowtypea`\arrowtypeb`\arrowtypec;\height>}}
\def\qtrianglep<#1>[#2`#3`#4;#5`#6`#7]{{
\settriparms[#1]
\width=\height                                
\diagram
\putqtrianglep<\arrowtypea`\arrowtypeb`
\arrowtypec;\height>
(0,0)[#2`#3`#4;#5`#6`{#7}]
\enddiagram
}}
 
\def\putdtrianglep<#1>(#2,#3)[#4`#5`#6;#7`#8`#9]{{%
\settriparms[#1]%
\xpos=#2 \ypos=#3
\puthmorphism(\xpos,\ypos)[#5`#6`{#9}]{\height}{\arrowtypec}b%
\advance\xpos by \height \advance\ypos by\height
\putmorphism(\xpos,\ypos)(-1,-1)[``{#7}]{\height}{\arrowtypea}l%
\putvmorphism(\xpos,\ypos)[#4``{#8}]{\height}{\arrowtypeb}r%
}}
 
\def\putdtriangle{\@ifnextchar <{\putdtrianglep}{\putdtrianglep
   <\arrowtypea`\arrowtypeb`\arrowtypec;\height>}}
\def\dtriangle{\@ifnextchar <{\dtrianglep}{\dtrianglep
   <\arrowtypea`\arrowtypeb`\arrowtypec;\height>}}
\def\dtrianglep<#1>[#2`#3`#4;#5`#6`#7]{{
\settriparms[#1]
\width=\height                                
\diagram
\putdtrianglep<\arrowtypea`\arrowtypeb`
\arrowtypec;\height>
(0,0)[#2`#3`#4;#5`#6`{#7}]
\enddiagram
}}
 
\def\putbtrianglep<#1>(#2,#3)[#4`#5`#6;#7`#8`#9]{{%
\settriparms[#1]%
\xpos=#2 \ypos=#3
\puthmorphism(\xpos,\ypos)[#5`#6`{#9}]{\height}{\arrowtypec}b%
\advance\ypos by\height
\putmorphism(\xpos,\ypos)(1,-1)[``{#8}]{\height}{\arrowtypeb}r%
\putvmorphism(\xpos,\ypos)[#4``{#7}]{\height}{\arrowtypea}l%
}}
 
\def\putbtriangle{\@ifnextchar <{\putbtrianglep}{\putbtrianglep
   <\arrowtypea`\arrowtypeb`\arrowtypec;\height>}}
\def\btriangle{\@ifnextchar <{\btrianglep}{\btrianglep
   <\arrowtypea`\arrowtypeb`\arrowtypec;\height>}}
\def\btrianglep<#1>[#2`#3`#4;#5`#6`#7]{{
\settriparms[#1]
\width=\height                               
\diagram
\putbtrianglep<\arrowtypea`\arrowtypeb`
\arrowtypec;\height>
(0,0)[#2`#3`#4;#5`#6`{#7}]
\enddiagram
}}
 
\def\putAtrianglep<#1>(#2,#3)[#4`#5`#6;#7`#8`#9]{{%
\settriparms[#1]%
\xpos=#2 \ypos=#3
{\multiply \height by2
\puthmorphism(\xpos,\ypos)[#5`#6`{#9}]{\height}{\arrowtypec}b}%
\advance\xpos by\height \advance\ypos by\height
\putmorphism(\xpos,\ypos)(-1,-1)[#4``{#7}]{\height}{\arrowtypea}l%
\putmorphism(\xpos,\ypos)(1,-1)[``{#8}]{\height}{\arrowtypeb}r%
}}
 
\def\putAtriangle{\@ifnextchar <{\putAtrianglep}{\putAtrianglep
   <\arrowtypea`\arrowtypeb`\arrowtypec;\height>}}
\def\Atriangle{\@ifnextchar <{\Atrianglep}{\Atrianglep
   <\arrowtypea`\arrowtypeb`\arrowtypec;\height>}}
\def\Atrianglep<#1>[#2`#3`#4;#5`#6`#7]{{
\settriparms[#1]
\width=\height                                     
\diagram
\putAtrianglep<\arrowtypea`\arrowtypeb`
\arrowtypec;\height>
(0,0)[#2`#3`#4;#5`#6`{#7}]
\enddiagram
}}
 
\def\putAtrianglepairp<#1>(#2)[#3;#4`#5`#6`#7`#8]{{%
\settripairparms[#1]%
\setpos(#2)%
\settokens`#3`%
\puthmorphism(\xpos,\ypos)[\tokenb`\tokenc`{#7}]{\height}{\arrowtyped}b%
\advance\xpos by\height
\puthmorphism(\xpos,\ypos)[\phantom{\tokenc}`\tokend`{#8}]%
{\height}{\arrowtypee}b%
\advance\ypos by\height
\putmorphism(\xpos,\ypos)(-1,-1)[\tokena``{#4}]{\height}{\arrowtypea}l%
\putvmorphism(\xpos,\ypos)[``{#5}]{\height}{\arrowtypeb}m%
\putmorphism(\xpos,\ypos)(1,-1)[``{#6}]{\height}{\arrowtypec}r%
}}
 
\def\putAtrianglepair{\@ifnextchar <{\putAtrianglepairp}{\putAtrianglepairp%
   <\arrowtypea`\arrowtypeb`\arrowtypec`\arrowtyped`\arrowtypee;\height>}}
\def\Atrianglepair{\@ifnextchar <{\Atrianglepairp}{\Atrianglepairp%
   <\arrowtypea`\arrowtypeb`\arrowtypec`\arrowtyped`\arrowtypee;\height>}}
 
\def\Atrianglepairp<#1>[#2;#3`#4`#5`#6`#7]{{
\settripairparms[#1]
\settokens`#2`
\width=\height                                
\diagram
\putAtrianglepairp                            
<\arrowtypea`\arrowtypeb`\arrowtypec`
\arrowtyped`\arrowtypee;\height>
(0,0)[{#2};#3`#4`#5`#6`{#7}]
\enddiagram
}}
 
\def\putVtrianglep<#1>(#2,#3)[#4`#5`#6;#7`#8`#9]{{%
\settriparms[#1]%
\xpos=#2 \ypos=#3
\advance\ypos by\height
{\multiply\height by2
\puthmorphism(\xpos,\ypos)[#4`#5`{#7}]{\height}{\arrowtypea}a}%
\putmorphism(\xpos,\ypos)(1,-1)[`#6`{#8}]{\height}{\arrowtypeb}l%
\advance\xpos by\height
\advance\xpos by\height
\putmorphism(\xpos,\ypos)(-1,-1)[``{#9}]{\height}{\arrowtypec}r%
}}
 
\def\putVtriangle{\@ifnextchar <{\putVtrianglep}{\putVtrianglep
   <\arrowtypea`\arrowtypeb`\arrowtypec;\height>}}
\def\Vtriangle{\@ifnextchar <{\Vtrianglep}{\Vtrianglep
   <\arrowtypea`\arrowtypeb`\arrowtypec;\height>}}
\def\Vtrianglep<#1>[#2`#3`#4;#5`#6`#7]{{
\settriparms[#1]
\width=\height                                 
\diagram
\putVtrianglep<\arrowtypea`\arrowtypeb`
\arrowtypec;\height>
(0,0)[#2`#3`#4;#5`#6`{#7}]
\enddiagram
}}
 
\def\putVtrianglepairp<#1>(#2)[#3;#4`#5`#6`#7`#8]{{
\settripairparms[#1]%
\setpos(#2)%
\settokens`#3`%
\advance\ypos by\height
\putmorphism(\xpos,\ypos)(1,-1)[`\tokend`{#6}]{\height}{\arrowtypec}l%
\puthmorphism(\xpos,\ypos)[\tokena`\tokenb`{#4}]{\height}{\arrowtypea}a%
\advance\xpos by\height
\puthmorphism(\xpos,\ypos)[\phantom{\tokenb}`\tokenc`{#5}]%
{\height}{\arrowtypeb}a%
\putvmorphism(\xpos,\ypos)[``{#7}]{\height}{\arrowtyped}m%
\advance\xpos by\height
\putmorphism(\xpos,\ypos)(-1,-1)[``{#8}]{\height}{\arrowtypee}r%
}}
 
\def\putVtrianglepair{\@ifnextchar <{\putVtrianglepairp}{\putVtrianglepairp%
    <\arrowtypea`\arrowtypeb`\arrowtypec`\arrowtyped`\arrowtypee;\height>}}
\def\Vtrianglepair{\@ifnextchar <{\Vtrianglepairp}{\Vtrianglepairp%
    <\arrowtypea`\arrowtypeb`\arrowtypec`\arrowtyped`\arrowtypee;\height>}}
\def\Vtrianglepairp<#1>[#2;#3`#4`#5`#6`#7]{{
\settripairparms[#1]
\settokens`#2`
\diagram
\putVtrianglepairp                             
<\arrowtypea`\arrowtypeb`\arrowtypec`
\arrowtyped`\arrowtypee;\height>
(0,0)[{#2};#3`#4`#5`#6`{#7}]
\enddiagram
}}

\def\putCtrianglep<#1>(#2,#3)[#4`#5`#6;#7`#8`#9]{{%
\settriparms[#1]%
\xpos=#2 \ypos=#3
\advance\ypos by\height
\putmorphism(\xpos,\ypos)(1,-1)[``{#9}]{\height}{\arrowtypec}l%
\advance\xpos by\height
\advance\ypos by\height
\putmorphism(\xpos,\ypos)(-1,-1)[#4`#5`{#7}]{\height}{\arrowtypea}l%
{\multiply\height by 2
\putvmorphism(\xpos,\ypos)[`#6`{#8}]{\height}{\arrowtypeb}r}%
}}
 
\def\putCtriangle{\@ifnextchar <{\putCtrianglep}{\putCtrianglep
    <\arrowtypea`\arrowtypeb`\arrowtypec;\height>}}
\def\Ctriangle{\@ifnextchar <{\Ctrianglep}{\Ctrianglep
    <\arrowtypea`\arrowtypeb`\arrowtypec;\height>}}
\def\Ctrianglep<#1>[#2`#3`#4;#5`#6`#7]{{
\settriparms[#1]
\width=\height                               
\diagram
\putCtrianglep<\arrowtypea`\arrowtypeb`
\arrowtypec;\height>
(0,0)[#2`#3`#4;#5`#6`{#7}]
\enddiagram
}}                                           
\def\putDtrianglep<#1>(#2,#3)[#4`#5`#6;#7`#8`#9]{{%
\settriparms[#1]%
\xpos=#2 \ypos=#3
\advance\xpos by\height \advance\ypos by\height
\putmorphism(\xpos,\ypos)(-1,-1)[``{#9}]{\height}{\arrowtypec}r%
\advance\xpos by-\height \advance\ypos by\height
\putmorphism(\xpos,\ypos)(1,-1)[`#5`{#8}]{\height}{\arrowtypeb}r%
{\multiply\height by 2
\putvmorphism(\xpos,\ypos)[#4`#6`{#7}]{\height}{\arrowtypea}l}%
}}
 
\def\putDtriangle{\@ifnextchar <{\putDtrianglep}{\putDtrianglep
    <\arrowtypea`\arrowtypeb`\arrowtypec;\height>}}
\def\Dtriangle{\@ifnextchar <{\Dtrianglep}{\Dtrianglep
   <\arrowtypea`\arrowtypeb`\arrowtypec;\height>}}
\def\Dtrianglep<#1>[#2`#3`#4;#5`#6`#7]{{
\settriparms[#1]
\width=\height                              
\diagram
\putDtrianglep<\arrowtypea`\arrowtypeb`
\arrowtypec;\height>
(0,0)[#2`#3`#4;#5`#6`{#7}]
\enddiagram
}}                                          
\def\setrecparms[#1`#2]{\width=#1 \height=#2}%
 
\def\recursep<#1`#2>[#3;#4`#5`#6`#7`#8]{{\m@th
\width=#1 \height=#2
\settokens`#3`
\settowidth{\tempdimen}{$\tokena$}
\ifdim\tempdimen=0pt
  \savebox{\tempboxa}{\hbox{$\tokenb$}}%
  \savebox{\tempboxb}{\hbox{$\tokend$}}%
  \savebox{\tempboxc}{\hbox{$#6$}}%
\else
  \savebox{\tempboxa}{\hbox{$\hbox{$\tokena$}\times\hbox{$\tokenb$}$}}%
  \savebox{\tempboxb}{\hbox{$\hbox{$\tokena$}\times\hbox{$\tokend$}$}}%
  \savebox{\tempboxc}{\hbox{$\hbox{$\tokena$}\times\hbox{$#6$}$}}%
\fi
\ypos=\height
\divide\ypos by 2
\xpos=\ypos
\advance\xpos by \width
\bfig
\putCtrianglep<-1`1`1;\ypos>(0,0)[`\tokenc`;#5`#6`{#7}]%
\puthmorphism(\ypos,0)[\tokend`\usebox{\tempboxb}`{#8}]{\width}{-1}b%
\puthmorphism(\ypos,\height)[\tokenb`\usebox{\tempboxa}`{#4}]{\width}{-1}a%
\advance\ypos by \width
\putvmorphism(\ypos,\height)[``\usebox{\tempboxc}]{\height}1r%
\efig
}}
 
\def\recurse{\@ifnextchar <{\recursep}{\recursep<\width`\height>}}
 
\def\puttwohmorphisms(#1,#2)[#3`#4;#5`#6]#7#8#9{{%
\puthmorphism(#1,#2)[#3`#4`]{#7}0a
\ypos=#2
\advance\ypos by 20
\puthmorphism(#1,\ypos)[\phantom{#3}`\phantom{#4}`#5]{#7}{#8}a
\advance\ypos by -40
\puthmorphism(#1,\ypos)[\phantom{#3}`\phantom{#4}`#6]{#7}{#9}b
}}
 
\def\puttwovmorphisms(#1,#2)[#3`#4;#5`#6]#7#8#9{{%
\putvmorphism(#1,#2)[#3`#4`]{#7}0a
\xpos=#1
\advance\xpos by -20
\putvmorphism(\xpos,#2)[\phantom{#3}`\phantom{#4}`#5]{#7}{#8}l
\advance\xpos by 40
\putvmorphism(\xpos,#2)[\phantom{#3}`\phantom{#4}`#6]{#7}{#9}r
}}
 
\def\puthcoequalizer(#1)[#2`#3`#4;#5`#6`#7]#8#9{{%
\setpos(#1)%
\puttwohmorphisms(\xpos,\ypos)[#2`#3;#5`#6]{#8}11%
\advance\xpos by #8
\puthmorphism(\xpos,\ypos)[\phantom{#3}`#4`#7]{#8}1{#9}
}}
 
\def\putvcoequalizer(#1)[#2`#3`#4;#5`#6`#7]#8#9{{%
\setpos(#1)%
\puttwovmorphisms(\xpos,\ypos)[#2`#3;#5`#6]{#8}11%
\advance\ypos by -#8
\putvmorphism(\xpos,\ypos)[\phantom{#3}`#4`#7]{#8}1{#9}
}}
 
\def\putthreehmorphisms(#1)[#2`#3;#4`#5`#6]#7(#8)#9{{%
\setpos(#1) \settypes(#8)
\if a#9 %
     \vertsize{\tempcounta}{#5}%
     \vertsize{\tempcountb}{#6}%
     \ifnum \tempcounta<\tempcountb \tempcounta=\tempcountb \fi
\else
     \vertsize{\tempcounta}{#4}%
     \vertsize{\tempcountb}{#5}%
     \ifnum \tempcounta<\tempcountb \tempcounta=\tempcountb \fi
\fi
\advance \tempcounta by 60
\puthmorphism(\xpos,\ypos)[#2`#3`#5]{#7}{\arrowtypeb}{#9}
\advance\ypos by \tempcounta
\puthmorphism(\xpos,\ypos)[\phantom{#2}`\phantom{#3}`#4]{#7}{\arrowtypea}{#9}
\advance\ypos by -\tempcounta \advance\ypos by -\tempcounta
\puthmorphism(\xpos,\ypos)[\phantom{#2}`\phantom{#3}`#6]{#7}{\arrowtypec}{#9}
}}
 
\def\setarrowtoks[#1`#2`#3`#4`#5`#6]{%
\def\toka{#1}
\def\tokb{#2}
\def\tokc{#3}
\def\tokd{#4}
\def\toke{#5}
\def\tokf{#6}
}
\def\hex{\@ifnextchar <{\hexp}{\hexp<1000`400>}}
\def\hexp<#1`#2>[#3`#4`#5`#6`#7`#8;#9]{%
\setarrowtoks[#9]
\yext=#2 \advance \yext by #2
\xext=#1 \advance\xext by \yext
\bfig
\putCtriangle<-1`0`1;#2>(0,0)[`#5`;\tokb``\tokd]
\xext=#1 \yext=#2 \advance \yext by #2
\putsquare<1`0`0`1;\xext`\yext>(#2,0)[#3`#4`#7`#8;\toka```\tokf]
\advance \xext by #2
\putDtriangle<0`1`-1;#2>(\xext,0)[`#6`;`\tokc`\toke]
\efig
}
\makeatother

\newtheorem{thm}{Theorem}[section]

\newtheorem{cor}[thm]{Corollary}
\newtheorem{lma}[thm]{Lemma}
\newtheorem{prop}[thm]{Proposition}

\theoremstyle{definition}
\newtheorem{defn}[thm]{Definition}
\newtheorem{eg}[thm]{Examples}

\theoremstyle{remark}

\newcommand{\F}{{\mathbb F}}
\renewcommand{\L}{{\mathcal L}}
\newcommand{\N}{{\mathbb N}}

\newcommand{\X}{{\mathcal X}}
\newcommand{\Y}{{\mathcal Y}}
\newcommand{\Z}{{\mathcal Z}}

\newif\ifnote 

\begin{document}

\title[Boolean Monomial Dynamical Systems]{Boolean 
Monomial Dynamical Systems} 
\author{Omar Col\'on-Reyes}

\author{Reinhard Laubenbacher}
\address{Virginia Bioinformatics Institute, Virginia Tech, 
Blacksburg, VA 24061-0477, USA} 

\author{Bodo Pareigis}
\address{Mathematisches Institut der Universit\"at 
M\"unchen, Germany} 

\begin{abstract}
An important problem in the theory of finite dynamical systems is to link
the structure of a system with its dynamics.  This paper contains such 
a link for a family of nonlinear systems over the field with two elements.
For systems that can be described by monomials (including Boolean
AND systems), one can obtain information about the limit cycle structure
from the structure of the monomials.
In particular, the paper contains a sufficient condition for a monomial
system to have only fixed points as limit cycles.  This condition depends
on the cycle structure of the dependency graph of the system and can
be verified in polynomial time.
\end{abstract}

\thanks{The research in this paper was supported in part by Contract Nr. 78761-SOL-0343 
from Los Alamos National Laboratory.}

\date{February 4, 2004}
\maketitle


 \section{Introduction} 
 Finite dynamical systems are time-discrete dynamical 
systems on finite state sets. Well-known examples include 
cellular automata and Boolean networks, which have found 
broad applications in engineering, computer science, and, 
more recently, computational biology.  (See, e.g., [K, AO, 
CS] for biological applications.) More general multi-state 
systems have been used in control theory [G, LB, M1, M2], 
the design and analysis of computer simulations [BR, BMR1, 
BMR2, LP2], and in computational biology as well [LS]. One 
underlying mathematical question that is common to many of 
these applications is how to analyze the dynamics of the 
models without actually enumerating all state transitions, 
since enumeration has exponential complexity in the number 
of model variables.  The present paper is a contribution 
toward an answer to this question. 
  
 For our purposes, a finite dynamical system is a function 
$f:X\to X$, where $X$ is a finite set [LP1]. Many 
applications assume $X$ to be of the form $k^n$, where 
$n\geq 1$ and $k$ is a finite field, often the field with 
two elements.  In this paper we restrict ourselves to the 
case $X=k^n=\F_2^n$.  The dynamics is generated by 
iteration of $f$, with the variables being updated 
simultaneously.      In this paper we present a family of 
nonlinear systems for which the above question can be 
answered, that is, for which one can obtain information 
about the dynamics from the structure of the function.  The 
answer is given in terms of properties of the dependency 
graph of the system, that is, the directed graph that 
represents the dependence of the coordinate functions of 
$f$ on the other coordinates. 

 We assume that $f:k^n \to  k^n$, $n\geq 1$, is a finite 
dynamical system. The dynamics of $f$ is encoded in its 
{\em state space} $\mathcal S(f)$, which is a directed 
graph defined as follows. The vertices of $\mathcal S(f)$ 
are the $2^n$ elements of $k^n$.  There is a directed edge 
$a\rightarrow b$ in $\mathcal S(f)$ if $f(a)=b$.  In 
particular, a directed edge from a vertex to itself is 
admissible. That is, $\mathcal S(f)$ encodes all state 
transitions of $f$, and has the property that every vertex 
has out-degree exactly equal to $1$.  Each connected graph 
component of $\mathcal S(f)$ consists of a directed cycle, 
a so-called {\em limit cycle}, with a directed tree 
attached to each vertex in the cycle, consisting of the so-
called {\em transients}. 

 Observe that $f$ can be described in terms of its 
coordinate functions $f_i:k^n \to  k$, that is, 
$f=(f_1,\ldots ,f_n)$.  It is well known that, if $k$ is 
any finite field with $q$ elements and $g:k^n \to  k$ is 
any set-theoretic function, then $g$ can be represented by 
a polynomial in $k[x_1,\ldots ,x_n]$ [LN, p. 369].  This 
polynomial can be chosen uniquely so that any variable in 
it appears to a degree less than the number of elements in 
$k$. That is, for any $g$ there is a unique $h\in 
k[x_1,\ldots ,x_n]/\langle x_i^q-x_i|i=1,\ldots ,n\rangle$, 
such that $g(a)=h(a)$ for all $a\in k^n$.   Consequently, 
any finite dynamical system over a finite field can be 
represented as a polynomial system.  This is the point of 
view we take in this paper.  For $q=2$ this implies that 
every function can be represented by a square-free 
polynomial. 

 Observe further that any polynomial function over $k=\F_2$ 
with square-free monomials can be represented as a Boolean 
function, with multiplication corresponding to the logical 
AND, addition to the logical XOR, and negation as addition 
of the constant term $1$.  This implies that any finite 
dynamical system over $\F_2$ can be realized as a Boolean 
network.  For a polynomial system $f=(f_1,\ldots ,f_n)$ the 
problem to be investigated then becomes one of drawing 
conclusions about the structure of the state space 
$\mathcal S(f)$ from the structure of the $f_i$.  We first 
present a brief survey of existing results. 

 In the case of a linear system over an arbitrary finite 
field this question has a complete answer [H]. Let $A$ be a 
matrix representation of a linear system $f:k^n \to  k^n$.  
(That is, the $f_i$ are linear polynomials without constant 
term.)  Then the number of limit cycles and their length, 
as well as the structure of the transients, can be 
determined from the factorization of the characteristic 
polynomial of the matrix $A$. The structure of the limit 
cycles had been determined earlier by Elspas [E], and for 
affine systems by Milligan and Wilson [MW]. For nonlinear 
systems a very elegant general approach was proposed in 
[C].  A general nonlinear system can be embedded in a 
linear one, to which the author then applied techniques 
like in [E] to obtain the number and length of the limit 
cycles.  The drawback of this method is that, if the 
nonlinear system has dimension $n$ and the field has $q$ 
elements, then the linear system has dimension $q^n$.  It 
is also very difficult to see directly the effect of the 
specific nonlinear functions on the state space structure. 

 Very few results exist for general classes of nonlinear 
systems.  Some facts about nonlinear one-dimensional 
cellular automata can be found in Wolfram's work [W].  For 
sequentially updated systems, there are a some results in 
[BMR1].  For instance, it is shown there (Proposition 5) 
that sequentially updated Boolean NOR systems do not have 
any fixed points. In the one-dimensional case this problem 
has been studied extensively over the $p$-adic numbers, 
that is, the dynamics of univariate polynomials over 
$\mathbf Q_p$, viewed as a dynamical system. In particular, 
the case of monomials has been considered [KN], which are 
the focus of the present paper. 

 In this paper we focus on the class of nonlinear systems 
over $k=\F_2$ described by special types of polynomials, 
namely monomials. That is, we consider systems $f=(f_i)$, 
so that each $f_i$ is a polynomial of the form 
$x_{i_1}x_{i_2}\cdots x_{i_r}$, or a constant equal to $0$ 
or $1$.  This class includes all Boolean networks made up 
of AND functions. Associated to a general polynomial system 
one can construct its {\em dependency graph} $\mathcal 
D(f)$, whose vertices $v_1,\ldots ,v_n$ correspond to the 
variables of the $f_i$.  There is a directed arrow 
$v_i\rightarrow v_j$ if $x_j$ appears in $f_i$.  A special 
definition needs to be made to account for constant 
polynomials. For Boolean monomial systems the dependency 
graph in fact allows the unambiguous reconstruction of the 
system. The main results of this paper show that in this 
case the cycle structure of the state space $\mathcal S(f)$ 
can be determined exclusively from the dependency graph 
$\mathcal N(f)$, that is, from the structure of the $f_i$.  
The key role is played by a numerical invariant associated 
to a strongly connected directed graph, that is, a graph in 
which there exists a walk (a directed path) between any two 
vertices.  For such a graph one can define its {\em loop 
number} as the minimum of the distances of two walks from 
some vertex to itself.  (The number is the same no matter 
which vertex is chosen.)  It turns out that the dependency 
graph of a monomial system can be decomposed into strongly 
connected components whose loop numbers determine the 
structure of the limit cycles. If the loop number of every 
strongly connected component is one, then the state space 
has only fixed points as limit cycles, that is, $f$ is a 
fixed point system. 

 \section{Preliminaries on Systems and Dependency Graphs} 

 Let $f = (f_1, \ldots, f_n):\F_2^n \to \F_2^n$ be a 
Boolean monomial parallel update system, where the monomial 
functions are of the form 
 $$f_i = \alpha_i {x_1}^{\varepsilon_{1i}} \ldots 
{x_n}^{\varepsilon_{ni}}$$ 
 with $\alpha_i \in \{0,1\}$ and $\varepsilon_{ji}  \in 
\{0,1\}$. If $\alpha_i = 0$ we set all $\varepsilon_{ji} = 
0$. Let $f^m := f\circ f \circ \ldots \circ f$ be the 
 $m$-fold composition of the map $f$ with itself. We write 
$f^m = ({f^m}_1, \ldots, {f^m}_n)$. By definition we have 
 $${f^m}_i = \alpha_i {{f^{m-1}}_1}^{\varepsilon_{1i}} 
\ldots {{f^{m-1}}_n}^{\varepsilon_{ni}}. $$ 

 \begin{defn} 
 With $f$ we associate a digraph $\X$, called {\em 
dependency graph}, with vertex set $\{a_1, \ldots, a_n, 
\varepsilon\}$. There is a directed edge from $a_i$ to 
$a_j$ if $\alpha_i = 1$ and $x_j$ is a factor in $f_i$ 
(i.e. $\varepsilon_{ji} = 1$).  There is a  directed edge 
from $a_i$ to $\varepsilon$ if $\alpha_i = 0$ (i.e. $f_i = 
0$). 
 \end{defn} 

 Observe that loops $a_i \to a_i$ are permitted. They occur 
if $f_i$ has the factor $x_i$. If there is an edge $a_i \to 
\varepsilon$ ($f_i = 0$), then there is no edge $a_i \to 
a_j$ for any $j$. It is straightforward to see that the 
monomial system $f$ is completely described by the 
dependency graph $\X$. 

 We give two simple examples: 
 \begin{eg} \label{examples} 
 \begin{enumerate} 
 \item The system $f$ given by the function $f = (x_3, x_1 
\cdot x_4, x_4, x_1)$ has the following dependency graph 
  $$\bfig 
 \putmorphism(0, 300)(1, 0)[a_1`a_3`]{300}1a 
 \putmorphism(0, 0)(1, 0)[a_2`a_4`]{300}1a 
 \putmorphism(0, 300)(1, -1)[``]{300}{-1}l 
 \putmorphism(0, 300)(0, -1)[``]{300}{-1}l 
 \putmorphism(300, 300)(0, -1)[``]{300}1l 
 \efig$$ 
 and the state space given in Figure 
\ref{examplestatespace}. 
 \begin{figure}[!htp] 
 \centering 

 \includegraphics[totalheight=4cm]{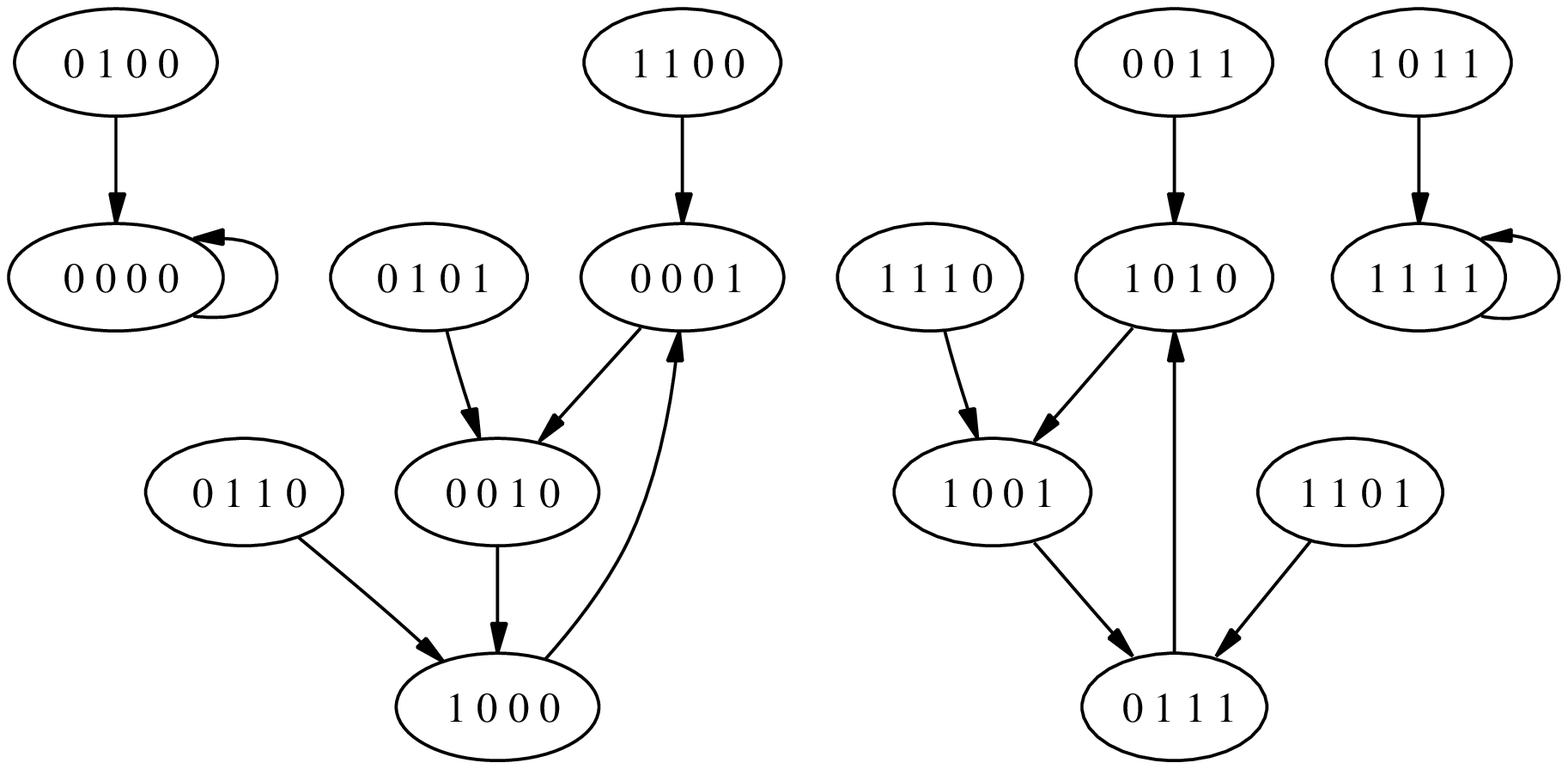} 
 \label{examplestatespace} 
 \caption{State space.} 
 \end{figure} 
 \item The system $f$ given by the function 
 $$f = (x_2, x_3 x_{10}, x_4x_7, x_5x_7, x_8x_9, x_3, 
x_6x_8, x_5x_9, x_8, x_8x_{11}, 0)$$ 
 has the following dependency graph 
  $$\bfig 
 \putmorphism(0, 600)(1, 0)[a_1`a_2`]{300}1a 
 \putmorphism(300, 600)(1, 0)[\phantom{a_2}`a_3`]{300}1a 
 \putmorphism(600, 600)(1, 0)[\phantom{a_3}`a_4`]{300}1a 
 \putmorphism(900, 600)(1, 0)[\phantom{a_4}`a_5`]{300}1a 
 \putmorphism(600, 300)(1, 0)[a_6`a_7`]{300}{-1}a 
 \putmorphism(900, 300)(1, 0)[\phantom{a_7}`a_8`]{300}1a 
 \putmorphism(1200, 320)(1, 0)[\phantom{a_8}`a_9`]{300}1a 
 \putmorphism(1200, 280)(1, 
0)[\phantom{a_8}`\phantom{a_9}`]{300}{-1}a 
 \putmorphism(300, 0)(1, 0)[a_{10}`a_{11}`]{300}1a 
 \putmorphism(600, 0)(1, 
0)[\phantom{a_{11}}`\varepsilon`]{300}1a 
 \putmorphism(0, 600)(1, -2)[``]{300}{-1}l 
 \putmorphism(600, 600)(1, -1)[``]{300}1l 
 \putmorphism(300, 600)(0, -1)[``]{600}1l 
 \putmorphism(600, 600)(0, -1)[``]{300}{-1}l 
 \putmorphism(900, 600)(0, -1)[``]{300}1l 
 \putmorphism(1180, 600)(0, -1)[``]{300}{-1}l 
 \putmorphism(1220, 600)(0, -1)[``]{300}1l 
 \putmorphism(1200, 600)(1, -1)[``]{300}1l 
 \putmorphism(1200, 300)(-3, -1)[``]{900}{-1}l 
 \efig$$ 
 Our theorems will show, that this system is a fixed point 
system, that is, all limit cycles are fixed points. In fact 
the state space of this system has 3 fixed points and 2048 
nodes. 
 \end{enumerate} 
 \end{eg} 

 \begin{prop} \label{xfactors} 
 Let $\X$ be the dependency graph of $f$ and assume 
${f^m}_i \not= 0$. There exists a walk $p: a_i \to a_j$ of 
length $m$ in $\X$ iff ${f^m}_i$ contains the factor $x_j$. 
 \end{prop} 

 \begin{proof} 
 Assume that there are directed edges from $a_i$ to 
$a_{j_1}, \ldots, a_{j_t}$. Any walk leaving $a_i$ goes 
through one of the $a_{j_1}, \ldots, a_{j_t}$ as a first 
step. So ${f^m}_i = \alpha_i {f^{m-1}}_{j_1} \ldots {f^{m-
1}}_{j_t}$. The claim now follows by induction and the 
observation that $x_j^2 = x_j$. For the converse observe 
that from 
 $$ {f^m}_i = \alpha_i {f^{m-1}}_{j_1} \ldots {f^{m-
1}}_{j_t}$$ 
 it follows that, if $x_j$ is a factor of ${f^m}_i$ it must 
also be a factor of some ${f^{m-1}}_{j_k}$. We can then 
again proceed by induction to get a walk of length $m$ from 
$a_i$ to $a_j$. 
 \end{proof} 

 \begin{cor} \label{functionproduct} 
 ${f^r}_i$ is the product of all functions ${f^s}_j$ for 
all walks $p:a_i \to a_j$ of length $s \leq r$. 
 \end{cor} 

 \begin{proof} 
 This follows by induction, as in Proposition 
\ref{xfactors}. 
 \end{proof} 

 \begin{cor} \label{zerocor} 
 If there is a walk $p: a_i \to a_j$ of length less than or 
equal to $r$, and if ${f^1}_j = f_j = 0$, then the function 
${f^r}_i$ is zero. 
 \end{cor} 

 \section{The structure of the dependency graph} 

 \begin{defn} 
 Let $\X$ be the dependency graph of $f$. 
 \begin{enumerate} 
 \item 
 We write $a \in \X$ for $a \in V_\X \setminus 
\{\varepsilon\}$. 
 \item 
 For vertices $a,b \in \X$, we call $a$ and $b$ {\em 
strongly connected}, and write $a \sim b$, if and only if 
there is a walk $p: a \to b$ and a walk $q: b \to a$. 
Observe that there is always a walk of length zero (the 
empty walk) from $a$ to $a$. Then $a \sim b$ is an 
equivalence relation on $V_\X \setminus \{\varepsilon\}$, 
called {\em strong equivalence}. 
 \item 
 The equivalence class of $a \in \X$ is called a {\em 
(strongly) connected component} and is denoted by $\bar a$. 
 Let $E(\X)$ be the set of equivalence classes. 
 \item 
 A vertex $a$ with an edge $a \to \varepsilon$ is called a 
{\em zero}. 
 \item 
 For $a, b \in \X$, let $p: a \to b$ be a walk. We denote 
the length of the walk $p$ by $|p|$. 
 \end{enumerate}  
 \end{defn} 

 The smallest strongly connected components are described 
as follows. If $a_i \in \X$ is a vertex with $f_i = 1$ then 
there is by definition no edge originating in $a_i$, so 
$a_i$ defines a one element strongly connected component 
which contains only the empty walk. The same holds in case 
$f_i = 0$, except that there is an edge $a_i \to 
\varepsilon$, but there is still only one walk, the empty 
walk, from $a_i$ to $a_i$. If $f_i = x_i$, then $a_i$ also 
defines a one element strongly connected component, since 
there is no edge $a_i \to a_j$ for $j \not= i$. However, 
there are infinitely many closed walks $p_j:a_i \to a_i$, 
one for each length $|p_j| = j \in \N$. 

 \begin{lma} \label{poset} 
 Define $\bar a \leq \bar b$ iff there is a walk from $a$ 
to $b$. Then $E(\X)$ is a poset. 
 \end{lma} 

 \begin{proof} 
 Observe that, given $a,a' \in \bar a$ and $b,b' \in \bar 
b$, there is a walk from $a$ to $b$ if and only if there is 
a walk from $a'$ to $b'$. If there is a walk from $a$ to 
$b$ and a walk from $b$ to $a$, then $a$ and $b$ are 
strongly connected and thus $\bar a = \bar b$. 
 \end{proof} 

 Observe that $a_i$ defines a one point connected component 
$\bar a_i$ if and only if $f_i = 0, 1,$ or $f_i=x_i$. 

 This partial order will come up again in the discussion of 
glueing. It covers all the edges that do not occur in 
strongly connected components. 

 \section{The state space of a connected component}  
\label{conncomp} 

 In this section we will study strongly connected 
components. The study of the relations between connected 
components, that lead to the poset structure of $E(\X)$ 
will be done in the next section. We first discuss three 
trivial cases. 

 Let $\X = \{a\}$ be a strongly connected component with 
$a$ a zero. Then $f^m = 0$ for all $m \geq 1$. Let $\X = 
\{a\}$ be a strongly connected component with $f = 1$. Then 
$f^m = 1$ for all $m \geq 1$. If $\X$ is the dependency 
graph of an arbitrary function $f: \F_2^r \to \F_2^r$ and 
if $a_i \in \X$ has a walk to a zero $a_j \to \varepsilon$ 
of length $n \geq 0$, then by Corollary \ref{zerocor} we 
have ${f^{n+m}}_i = 0$ for all $m \geq 1$. 

 Let $\bar a$ be the strongly connected component of $a \in 
\X$. Assume there is a walk from $a$ to a zero $a_j$ in 
$\X$. Then there is an $n \geq 0$ such that ${f^{n+m}}_i = 
0$ for all $a_i \in \bar a$ and all $m \geq 1$. In 
particular if $\X$ is strongly connected and if $a_j \in 
\X$ is a zero, then $f$ is a fixed point system with 
exactly one fixed point $(0, \ldots, 0)$. So for the rest 
of section \ref{conncomp} we will assume that $\X$ is 
strongly connected and we exclude the cases $f(x_1) = 1$ 
and $f(x_1) = 0$. 

 \subsection{Loop numbers} \label{loopnumbers} 

 \begin{defn} 
 The {\em loop number} of $a \in \X$ is the minimum of all 
numbers $t \geq 1$ with $t = |p| - |q|$ for all closed 
walks $p,q: a \to a$. If there is no closed walk from $a$ 
to $a$ then we set the loop number to be zero. This last 
case occurs only if $\bar a = \{a\}$ and there is no edge 
from $a$ to $a$ (i.e. $f = 0,1$, but we have excluded $f = 
0$.). 
 \end{defn} 

 If there is a loop $p: a \to a$ then the loop number of 
$a$ is $|pp| - |p| = 1$. 

 Let $a,b \in \X$. We show that the loop numbers of $a$ and 
of $b$ are equal. We get 

 \begin{lma} 
 The loop number is constant on any strongly connected 
$\X$, so the loop number of a strongly connected $\X$ is a 
well defined number $\L(\X)$. 
 \end{lma} 

 \begin{proof} 
 Let $p':a \to b$ and $q':b \to a$ be walks. Then $p'pq', 
p'qq': b \to b$ are closed walks with $|p'pq'|- |p'qq'| = 
t$, so the loop number of $b$ is less than or equal to the 
loop number of $a$. By symmetry the loop number is constant 
on $\X$. 
 \end{proof} 

 \begin{lma} \label{L4.3} 
 Let the loop number of $\X$ be $t$. Let $p': a_i \to a_j$ 
and $q': a_i \to a_j$ be walks. Then $|p'| - |q'| \in (t) 
\subseteq {\mathbb Z}$. 
 \end{lma} 

 \begin{proof} 
 Assume $|p'| > |q'|$ and let $|p'| - |q'| = rt + s$ with 
$0 \leq s < t$. We want to show $s=0$. Let $p,q : a_i \to 
a_i$ be such that $|q| - |p| = t$. We have $r\geq 0$. Then 
 $$|p'p| - |q'q| = |p'| + |p| - |q'| - |q| = rt + s -t = 
(r-1)t + s.$$ 
 Hence there are walks $p'',q'': a_i \to a_j$ with $|p''| - 
|q''| = s$. Let $p^*: a_j \to a_i$ be a walk. Then 
$|p^*p''| - |p^*q''| = s = 0$ because of the minimality of 
the loop number $t$. So $|p'| - |q'| \in (t)$. 
 \end{proof} 

 \begin{cor} \label{looplength} 
 Let the loop number of $\X$ be $t$. Let $p:a \to a$ be a 
closed walk. Then $|p| \in (t)$. 
 \end{cor} 

 \begin{proof} 
 Im Lemma \ref{L4.3} take $p$ and $pp$. 
 \end{proof} 

 \begin{prop} \label{pathstab} 
 Let the loop number of $\X$ be $t \geq 1$. For each $a,b 
\in \X$ there exists an $m \in \N$ and walks $p_i:a \to b$ 
of length $|p_i| = m+i\cdot t$ for all $i \in \N$.  
 \end{prop} 

 \begin{proof} 
 Let $p,q:a \to a$ be closed walks with $|q|-|p|=t$. By 
Corollary \ref{looplength} $|p|$ is divisible by $t$. Let 
$r := |p|/t$. Let $p':a \to b$ be a walk of length $s := 
|p'|$. Let $m := s + (r^2 - r)t$. Write $i \in \N$ as $i = 
jr + k$ with $0 \leq k < r$. Then the composition of walks 
 $$p'q^kp^{r-1+j-k}: a \to b$$ 
 has length $|p'q^kp^{r-1+j-k}| 
 \ifnote \footnote{$ 
 = s + k(r+1)t + (r-1+j-k)rt = s + (kr+k+r^2-r+jr-kr)t = s 
+ (r^2 - r)t + (jr + k)t $} 
 \fi 
 = m + i\cdot t$. 
 \end{proof} 

 Observe that the number $m$ can be quite large. Even if we 
take $p'$ to be a walk of minimal length and $p,q$ of 
minimal length such that $|q|-|p| = t$, i.e. $s$ and $r$ 
minimal, we get $m = s + (r^2 - r)t$. This contributes to 
the lengths of transients in the state space. 

 \subsection{Fixed points and cycles of strongly connected 
components} 

 \begin{lma} 
 For $a_i, a_j \in \X$ we define 
 $$a_i \approx a_j :\iff \exists \mbox{ walk } p: a_i \to 
a_j \mbox{ with } |p| \in (t).$$ 
 This is an equivalence relation, called {\em loop 
equivalence}. 
 \end{lma} 

 \begin{proof} 
 We have $a_i \approx a_i$ with a walk of length zero and 
thus reflexivity. Transitivity is trivial. For symmetry let 
$a_i \approx a_j$, then there is a walk $p:  a_i \to a_j$ 
with $|p| \in (t)$. Let $q: a_j \to a_i$ be any walk. Then 
$qp: a_i \to a_i$ is a walk with $|qp| = |q| + |p| \in (t)$ 
by Corollary \ref{looplength}, hence $|q| \in (t)$ and thus 
$a_j \approx a_i$. 
 \end{proof} 

 \begin{lma} 
 There are exactly $t$ loop equivalence classes in $\X$. 
 \end{lma} 

 \begin{proof} 
 Since there are walks starting in $a_i$ for all lengths 
$\geq 0$, take a walk $a_i \to a_{i+1} \to \ldots \to 
a_{i+t}$ of length $t$. The $a_{i+j}$, $j=0, \ldots, t-1$ 
are in different equivalence classes 
 $$\bar{\bar{a}}_i, \bar{\bar{a}}_{i+1}, \ldots, 
\bar{\bar{a}}_{i+t} = \bar{\bar{a}}_i$$ 
 by Lemma \ref{L4.3}. Every $a_k$ is in one of these 
equivalence classes, since there is a walk $q:a_{i+t} \to 
a_k$ of length $|q| = rt + s$, so there is a walk $q': 
a_{i+s} \to a_k$ of length $|q'| \in (t)$, hence $a_{i+s} 
\approx a_k$. 
 \end{proof} 

 We may now label and enumerate the vertices of $\X$ in the 
following way    
 $$(a_1, \ldots, a_{i_1}),(a_{i_1+1}, \ldots, a_{i_2}), 
\ldots, (a_{i_{t-1}+1}, \ldots, a_{i_t})$$ 
 where each group $(a_{i_j+1}, \ldots, a_{i_{j+1}})$ is a 
loop equivalence class and there is an edge $a_{i_j} \to 
a_{i_j+1}$ for all $j = 1, \ldots, t-1$. If $s_j := i_j - 
i_{j-1}$ is the number of elements in each loop class, then 
$\sum_{i=1}^t s_i = n$. The following example of a strongly 
connected dependency graph of loop number 3 should  
illuminate this enumeration of vertices 
$(a_1,a_2,a_3),(a_4,a_5),(a_6,a_7)$: 
  $$\bfig 
 \putmorphism(0, 600)(1, 0)[a_1`a_7`]{300}{-1}a 
 \putmorphism(300, 600)(1, 0)[\phantom{a_2}`a_2`]{300}1a 
 \putmorphism(0, 600)(1, -2)[`a_4`]{150}1l 
 \putmorphism(300, 600)(-1, -2)[``]{150}{-1}l 
 \putmorphism(600, 600)(1, -2)[`a_5`]{150}1l 
 \putmorphism(300, 0)(1, 0)[a_3`a_6`]{300}{-1}a 
 \putmorphism(750, 300)(-1, -2)[``]{150}1l 
 \putmorphism(150, 300)(1, -2)[``]{150}{-1}l 
 \efig$$ 

 \begin{cor} \label{cycleclassm} 
 Let $a \in \X$ and let $\{a_1, \ldots, a_u\}$ be the loop 
class of $a$ (with $a_1 = a$). Then there is an $m \in \N$ 
such that for all $j = 1, \ldots, u$ and all $i \in \N$ 
there is a walk $a \to a_j$ of length $(m+i)t$. 
 \end{cor} 

 \begin{proof} 
 Use Proposition \ref{pathstab} to obtain $m_j$ and walks 
$p_{ij}: a\to a_j$ of lengths $m_j + it$ for all $i \in \N$ 
and all $j = 1, \ldots, u$. By Lemma \ref{L4.3} each $m_j$ 
is divisible by $t$. Take $m := \max\{m_1, \ldots, 
m_u\}/t$. 
 \end{proof} 

 Before we study the fixed points and cycles of an 
arbitrary connected dependency graph we look at a simple 
example. 
 
 \begin{prop} 
 The state space of a directed $t$-gon is isomorphic to the 
set of orbits of the action of the cyclic group of order 
$t$ acting on $\F_2^t$, the $t$-dimensional hypercube, by 
cyclically exchanging the canonical basis vectors. 
 \end{prop} 

 \begin{proof} 
 The dynamical system $f$ of a directed $t$-gon is $f(x_1, 
\ldots,x_t) = (x_2, x_3, \ldots, x_1)$. The states in the 
state space are elements of $\F_2^t$, and the action of $f$ 
and the powers of $f$ give an action of the cyclic group 
$C_t$ on $\F_2^t$ by cyclic exchange of the basis vectors. 
 \end{proof} 

 \begin{thm} \label{main1} 
 Let $\X$ be strongly connected with loop number $t \geq 1$ 
and $n$ vertices. Let    
 $$(a_1, \ldots, a_{i_1}),(a_{i_1+1}, \ldots, a_{i_2}), 
\ldots, (a_{i_{t-1}+1}, \ldots, a_{i_t})$$ 
 be the enumeration of vertices of $\X$ as described above. 
Then there is an $m$ such that 
 $$f^{mt} = (\underbrace{y_1, \ldots, y_1}_{s_1{\mbox{ 
\tiny times}}}, \underbrace{y_2, \ldots, y_2}_{s_2\mbox{ 
\tiny times}}, \ldots, \underbrace{y_t, \ldots, 
y_t}_{s_t\mbox{ \tiny times}})$$ 
 where $y_1 = x_1 \cdot\ldots\cdot x_{i_1}$, $y_2 = 
x_{i_1+1}\cdot \ldots \cdot x_{i_2}$, $\ldots$, $y_t = 
x_{i_{t-1}+1}\cdot \ldots \cdot x_{i_t}$. 

 Furthermore 
 $$f^{mt+1} = (\underbrace{y_2, \ldots, y_2}_{s_1{\mbox{ 
\tiny times}}}, \underbrace{y_3, \ldots, y_3}_{s_2{ 
\mbox{\tiny times}}}, \ldots, \underbrace{y_1, \ldots, 
y_1}_{s_t{\mbox{ \tiny times}}}).$$ 
 $$\vdots$$ 
 $$f^{mt+j} = (\underbrace{y_{j+1}, \ldots, 
y_{j+1}}_{s_1{\mbox{ \tiny times}}}, \underbrace{y_{j+2}, 
\ldots, y_{j+2}}_{s_2{ \mbox{\tiny times}}}, \ldots, 
\underbrace{y_{j}, \ldots, y_{t}}_{s_t{\mbox{ \tiny 
times}}}).$$ 
 \end{thm} 

 \begin{proof} 
 For each loop class $c$ we use some $m_c$ from Corollary 
\ref{cycleclassm} and compute the values ${f^{mt}}_i$ for 
all $i$. By Proposition \ref{xfactors} we get ${f^{mt}}_i = 
x_1 \ldots x_u$ for all $i = 1, \ldots, u$. If we use the 
maximum $m$ of all such $m_c$, then we have the structure 
of $f^{mt}$. 

 To determine ${f^{mt+1}}_i$ observe that for $i = 1, 
\ldots, i_1$ there are walks $a_i \to a_j$ of length $mt+1$ 
for all $a_j \in \{a_{i_1+1}, \ldots, a_{i_2}\}$. By Lemma 
\ref{L4.3} no walk of length $mt+1$ starting in $a_i$ can 
end in a vertex different from these $a_j$. Hence 
${f^{mt+1}}_i = y_2$ for all $i = 1, \ldots, i_1$. The rest 
of the proof follows by induction. 
 \end{proof} 

 This theorem allows us to give a complete description of 
the cycles in the state space of $\X$. 

 \begin{cor} 
 Let $\X$ be as in Theorem \ref{main1}. Then the subgraph 
of cycles in the state space of $f$ is isomorphic to the 
state space of a directed $t$-gon, hence the set of orbits 
in the hypercube $\F_2^t$ under the action of the cyclic 
group $C_t$. 
 \end{cor} 

 \begin{proof} 
 For every choice of arguments for the $x_i$ we end up in 
the form of $f^{mt}$. Then $f$ acts on these points by a 
cyclic permutation, which defines the various cycles in the 
state space. By reducing the number of $y$'s with the same 
subscript to one, we get the system $(y_1, \ldots, y_t)$ 
with the same cyclic action $g(y_1, \ldots, y_n) = (y_2, 
y_3, \ldots, y_1)$ and $g^j(y_1, \ldots, y_n) = (y_{j+1}, 
y_{j+2}, \ldots, y_j)$. This system arises from a directed 
$t$-gon. 
 \end{proof} 

 \begin{cor} \label{fixedpoint1} 
 Let $\X$ be as in the Theorem. 
 \begin{enumerate} 
 \item The system $f$ has fixed points $(0, \ldots, 0)$ and 
$(1, \ldots, 1)$. 
 \item The system $f$ has limit cycles of all lengths dividing $t$.  
 \item The system $f$ is a fixed point system if and only 
if the loop number of $\X$ is 1. 
 \end{enumerate} 
 \end{cor} 

 Example \ref{examples} (1) has two strongly connected 
components $\{a_2\}$ and $\{a_1, a_3, a_4\}$ with loop 
numbers 0 and 3 respectively. The second strongly connected 
component has two 3-cycles as well as two fixed points. 
 
 Example \ref{examples} (2) has four strongly connected 
components $\{a_1, a_2, a_{10}\}$, $\{a_3, a_4, a_6, 
a_7\}$,\break $\{a_5, a_8, a_9\}$, and $\{a_{11}\}$ with 
loop numbers 3, 1, 1, and 0 respectively. The second 
strongly connected component has two 3 cycles and two fixed 
points. The next two components have two fixed points each, 
and the last component has one fixed point. In both 
examples the strongly connected components of the systems 
are connected by further edges. So we have to investigate, 
how these additional edges effect the fixed point property. 

 \subsection{The computation of loop numbers} 

 In view of the importance of the loop number for 
determining if a system is a fixed point system, we 
describe some polynomial time algorithm to compute the loop 
number of a strongly connected component. Let $A$ be the 
adjacency matrix of an arbitrary dependency graph $\X$ with 
$n$ vertices. Then the connected components can be read off 
the powers $A, A^2, \ldots, A^n$ of $A$ by the fact that a 
vertex $a_i$ is connected by a walk of length $s$ to a 
vertex $a_j$ iff the $ij$-th component of $A^s$ is non-
zero. 

 The strongly connected component of a vertex $a_i$ is 
obtained as follows: take the $i$-th row of all the 
matrices $A^r$, find all nonzero entries and determine the 
set $R(i)$ of associated $j$s. Take the $i$-th column of 
all the matrices $A^r$, find all nonzero entries and 
determine the set $C(i)$ of associated $j$s. Then $R(i) 
\cap C(i)$ is the set of indices $j$ such that $a_j$ is 
strongly connected with $a_i$. 

 Now we assume that $\X$ is strongly connected and we 
exclude the cases $f(x_1) = 1$ and $f(x_1) = 0$. A closed 
walk $p: a_i \to a_i$ is called a {\em circuit} in $\X$ if 
it has no repetitive vertices. A circuit has length at most 
$n$, so it can be read off the powers $A, A^2, \ldots, A^n$ 
of the adjacency matrix $A$. In contrast, closed walks may 
have arbitrarily large lengths. 

 \begin{thm} 
 Let $X$ be a strongly connected graph. The loop number of 
$\X$ is the greatest common divisor of the numbers $i$ with 
$1 \leq i \leq n$, such that $A^i$ has at least one non-
zero diagonal entry. 
 \end{thm} 

 \begin{proof} 
 We first prove that the loop number of $\X$ is the 
greatest common divisor of the lengths of all circuits of 
$\X$. The problem in proving this is, that the loop number 
cannot be represented in general as the difference of the 
lengths of two circuits. 

 Observe that a closed walk $p$ can be decomposed into a 
number of circuits $p_1, \ldots, p_l$ sharing vertices. 

 Let $d$ be the greatest common divisor of the lengths of 
all circuits of $\X$. Take a vertex $a \in \X$. Let $p$ and 
$q$ be closed walks through $a$ representing the loop 
number $|p| - |q| = t$ of $\X$. We want to show $d = t$. 

 Decompose $p$ and $q$ into a number of circuits $p_1, 
\ldots, p_l$ and $q_1, \ldots, q_m$. Then we get $|p| = 
|p_1| + \ldots +  |p_l|$ and $|q| = |q_1| + \ldots +  
|q_m|$. Hence 
 $$|p_1| + \ldots +  |p_l| - |q_1| - \ldots - |q_m| = t,$$ 
 so that $d$ divides $t$. 

 Now we show that $t$ divides the length of each circuit in 
$\X$. Assume there is a circuit $p'$ through $b \in \X$ 
whose length $s := |p'|$ is not divisible by $t$. There are 
walks $q_1: a \to b$ and $q_2:b \to a$. Let $r := |q_1| + 
|q_2|$ and let $t' := \gcd(r,s,t)$. Then $t'$ can be 
written as $t' = \alpha t - \beta r - \gamma s$ with 
$\alpha, \gamma \geq 0$ and $\beta > 0$. So 
 $$|p^\alpha| - |q^\alpha q_2 (q_1q_2)^{\beta-1} 
(p')^\gamma q_1| = t',$$ 
 i.e., we have constructed two closed walks whose lengths 
differ by $t'$.  Hence $t' = t$ is the loop number of $\X$. 
Then $t' = t$ divides $s = |p'|$, a contradiction. Thus all 
lengths of circuits in $\X$ are divisible by $t$, so that 
$t$ divides $d$, hence $d = t$. 

 Now the length of all possible closed walks (of lengths 
$\leq n$) can be read of the non-zero diagonal entries of 
the powers of the adjacency matrix. These walks are 
composed, as above, of circuits, and their lengths are sums 
of the lengths of certain circuits. So the greatest common 
divisor of the numbers $i$, such that $A^i$ has at least 
one non-zero diagonal entry is the greatest common divisor 
of the lengths of all circuits of $\X$. 
 \end{proof} 
 
 \section{Glueing} 

 \begin{defn} 
 Let $\X$ and $\Y$ be dependency graphs of functions $f: 
\F_2^r \to \F_2^r$ and $g: \F_2^s \to \F_2^s$, 
respectively. A {\em glueing $\X \# \Y$ of $\Y$ to $\X$} 
consists of a digraph with vertices $V_\X \dot\cup V_\Y$ 
and edges $E_\X \dot\cup E_\Y$ (disjoint union), together 
with a set of additional directed edges from vertices in 
$\Y$ to vertices in $\X$. 

 The function of $\X \# \Y$ is denoted by $f \# g: 
\F_2^{r+s} \to \F_2^{r+s}$. 
 \end{defn} 

 Observe that $\X$ is a subgraph of $\X \# \Y$ and that $f: 
\F_2^r \to \F_2^r$ is a quotient system of $f\#g: 
\F_2^{r+s} \to \F_2^{r+s}$. 

 We write elements in $\F_2^{r+s} = \F_2^r \times \F_2^s$ 
as pairs $(\alpha,\beta)$ with $\alpha = (\alpha_1, \ldots, 
\alpha_r)$ and $\beta = (\beta_1, \ldots, \beta_s)$. 
Similarly we write the variables for the function $f\#g$ as 
$(x,y)$ with $x = (x_1, \ldots,x_r)$ and $y = (y_1, \ldots, 
y_s)$. Thus we can write the function $f\#g$ as 
 $$(f\#g)(x_1, \ldots, x_r, y_1, \ldots, y_s) = (f(x_1, 
\ldots, x_r), h(x_1, \ldots, x_r, y_1, \ldots, y_s)),$$ 
 with $f = (f_1, \ldots, f_r)$ and $h = (h_1, \ldots, h_s) 
= ((f\#g)_{r+1}, \ldots, (f\#g)_{r+s})$. 

 We will use the poset $E(\X)$ as in Lemma \ref{poset} 
together with the glueing procedure to study fixed point 
systems. We will discuss elements in $E(\X)$ with no edges 
($f = 0,1$) separately. They are the strongly connected 
components of loop length 0. 

 \begin{lma} \label{fpscomponents} 
 Let $\X \# \Y$ be a fixed point system. Then $\X$ is a 
fixed point system. 
 \end{lma} 

 \begin{proof} 
 This follows from the fact that $f$ is a quotient system 
of $f \# g$. Indeed let $(\alpha,\beta) \in \F_2^{r+s}$. 
Let $m \in \N$ be such that $(f \# g)^m(\alpha,\beta) = (f 
\# g)^{m+1}(\alpha,\beta)$. Then we have $f^m(\alpha) = 
f^{m+1}(\alpha)$ so that $f$ is a fixed point system. 
 \end{proof} 

 \begin{thm} \label{main2} 
 Let $\X$ be a fixed point system and let $\Y$ be strongly 
connected of loop length $\geq 1$.  Let $\X \# \Y$ be a 
glueing. 
 The following are equivalent: 
 \begin{itemize} 
 \item The glueing $\X \# \Y$ is a fixed point system. 
 \item 
  \begin{enumerate} 
  \item There is a vertex $a\in \Y$ that is connected with 
a walk to a zero in $\X$, or 
  \item $\Y$ is a fixed point system. 
  \end{enumerate} 
 \end{itemize} 
 \end{thm} 

 \begin{proof} 
 $\Longleftarrow:$ We use the notation introduced above. 
Assume that there is a vertex $a\in \Y$ that is connected 
with a walk to a zero in $\X$. Then all vertices of $\Y$ 
are connected to a zero. By Corollary \ref{zerocor} we get 
an $n \in \N$ such that ${(f\#g)^{n+m}}_{r+i} = 0$ for all 
$m \geq 1$ and for all $i = 1, \ldots, s$. Since $f$ is a 
fixed point system we have $f^m = f^{m+1} = \ldots$. So for 
a sufficiently large exponent $m'$ the component $\Y$ in 
$\X \# \Y$ can only contribute a fixed point $\beta = (0, 
\ldots, 0)$. Hence $(f \# g)^{m'} = (f \# g)^{m'+1}$. 

 Let $g$ be a fixed point system. Note that $f$ is a 
subsystem of $f \# g$. Iterate $f \# g$ until the component 
on $\X$ becomes constant: 
 $$(f \# g)^m(x,y) = (f^m(x),({[(f \# g)^m]}_{r+1}(x,y), 
\ldots, {[(f \# g)^m]}_{r+s}(x,y)),$$ 
 where 
 $$(f \# g)^{m+k}(x,y) = (f^m(x),({(f \# 
g)^{m+k}}_{r+1}(x,y), \ldots, {(f \# 
g)^{m+k}}_{r+s}(x,y)).$$ 
 Set 
 $$(z_1, \ldots, z_s) := ({[(f \# g)^m]}_{r+1}(x,y), 
\ldots, {[(f \# g)^m]}_{r+s}(x,y)).$$ 
 Then apply powers of $f \# g$ to $(f^m(x),(z_1, \ldots, 
z_s))$ to get 
 $$(f \# g)^{m+1}(x,y) = (f^m(x),(f \# g)_{r+1}(f^m(x),z), 
\ldots, (f \# g)_{r+s}(f^m(x),z)),$$ 
 and 
 $$(f \# g)^{m+k}(x,y) = (f^m(x),{(f \# 
g)^k}_{r+1}(f^m(x),z), \ldots, {(f \# 
g)^k}_{r+s}(f^m(x),z)).$$ 
 So, for a fixed choice of $x$, the functions $(f \# 
g)_{r+1}, \ldots, (f \# g)_{r+s}$, when viewed only as a 
function on elements from $\F_2^s$, are the functions $g_1, 
\ldots, g_s$ multiplied with certain factors from the fixed 
point $f^m(x)$. This defines a new system $h_x: \F_2^s \to 
\F_2^s$, so that 
 $$(f \# g)^k(f^m(x),z) = (f^m(x), h_x^k(z)).\eqno{(*)}$$ 
 If one of the factors taken from $f^m(x)$ is zero, then by 
Corollary \ref{zerocor} and by the fact that $\Y$ is 
strongly connected, we get that $h_x$ is a fixed point 
system with the only fixed point $\beta = (0, \ldots, 0)$. 
If all of the factors taken from $f^m(x)$ are 1, then $h_x 
= g$ and hence is a fixed point system. So by (*) we see 
that $(f \# g)^k(f^m(x),z)$ ends in a fixed point for all 
choices of $x$ and $y$. 

 $ \Longrightarrow :$ Before we give the proof in this 
direction, we prove a Lemma. Let $\Z$ be the dependency 
graph of an arbitrary Boolean monomial system. Decompose 
$\Z$ into two parts, where $\Z_0$ is the glueing of all 
strongly connected components that allow a walk to a zero. 
Let $\Z_1$ be the glueing of all the other strongly 
connected components of $\Z$. Then there are no walks from 
any vertex in $\Z_1$ to any vertex in $\Z_0$, so that $\Z$ 
is a glueing of $\Z_0$ to $\Z_1$. 

 \begin{lma} 
 Any system of the form $\Z_1$ has a fixed point $(1, 
\ldots, 1)$. 
 \end{lma} 

 \begin{proof} 
 This is proved by induction on the number of connected 
components of $\Z_1$. Assume that $\Z_1 = \X_1 \# \Y$ where 
$\Y$ is a connected component with no zero. We can assume 
that $\X_1$ has a fixed point $(1, \ldots, 1)$. Observe 
that $\Y$ belongs to the component $(\X \# \Y)_1$ (not 
connected to a zero) and $\Y$ has a fixed point $(1, 
\ldots, 1)$ by Corollary \ref{fixedpoint1}. So we get from 
the induction hypothesis that $\Z_1$ has also a fixed point 
$(1, \ldots, 1)$. 
 \end{proof} 

 {\em Proof of Theorem continued.} Now assume that $\X \# 
\Y$ is a fixed point system, and that no vertex of $\Y$ is 
connected with a walk to a zero in $\X$. Now we use as  
initial state $(x,y)$ for the system $f \# g$, with $x = 
(1, \ldots, 1, 0, \ldots, 0)$ (a fixed point of $f$), where 
the component $(1, \ldots,1)$ belongs to $\X_1$ and $(0, 
\ldots, 0)$ belongs to $\X_0$, and $y$ arbitrary.  Then, as 
discussed above, the system $h_x: \F_2^s \to \F_2^s$ 
coincides with $g: \F_2^s \to \F_2^s$. Since $f \# g$ is a 
fixed point system, the system $g$ can also only contain 
fixed points and no proper limit cycles. Thus $g$ is a 
fixed point system. 
 \end{proof} 

 \begin{thm} \label{main3} 
 Let $f$ be a fixed point system and let $\Y$ be strongly 
connected of loop length 0.  Let $\X \# \Y$ be a glueing. 
Then the system $f\# g$ corresponding to the glueing $\X \# 
\Y$ is a fixed point system. 
 \end{thm} 

 \begin{proof} 
 There are two cases: $\Y$ is a one point component with 
$f_a = 1$ and the case that $\Y$ is a one point component 
with $f_a = 0$. If $f_a = 0$, then there is no additional 
edge from $\Y = \{a\}$ to $\X$ and thus $f \# g$ is a fixed 
point system with fixed points of the form $(\alpha,0)$, 
where $\alpha$ is a fixed point of $f$. If $f_a = 1$ then 
any additional edge from $\Y = \{a\}$ to $a_i \in \X$ adds 
a factor to the last component, so that $f \# g = (f_1, 
\ldots, f_n, x_{i_1}\ldots x_{i_r})$. If $\alpha = 
(\alpha_1, \ldots, \alpha_n)$ is a fixed point for $f$, 
i.e. $f(\alpha) = \alpha)$ then $(f \# g)(\alpha,\beta) = 
(\alpha, \alpha_{i_1}\cdots \alpha_{i_r})$ which is again a 
fixed point. So $f\#g$ is a fixed point system. 
 \end{proof} 

 \section{Boolean monomial fixed point systems} 

 In this section $\X$ will be the digraph of an arbitrary 
Boolean monomial parallel update system. 

 \begin{thm} \label{main4} 
 Let $\X$ be the digraph of $f: \F_2^n \to \F_2^n$. The 
following are equivalent: 
 \begin{enumerate} 
 \item The system $f: \F_2^n \to \F_2^n$ is a fixed point 
system. 
 \item For every vertex $a \in \X$ one of the following 
holds 
  \begin{enumerate} 
  \item $a$ allows two closed walks $p,q: a \to a$ of 
length $|q| = |p| + 1$, 
  \item $a$ is connected with a walk to a zero, or 
  \item there is no walk of length $\geq 1$ from $a$ to 
$a$. 
  \end{enumerate} 
 \end{enumerate} 
 \end{thm} 

 \begin{proof} 
 This is an immediate consequence of Theorem \ref{main2}, 
Corollary \ref{fixedpoint1} and the remarks at the 
beginning of section \ref{conncomp}. 
 \end{proof} 

 \begin{defn} A system 
 $$ f=(f_{1},f_{2},...,f_{n}):\F_{2}^{n}\rightarrow 
\F_{2}^{n} $$ 
 is a {\em triangular system}, if each $f_i$ is of the form 
 $$ f_{i}=\alpha_{i}x_{1}^{\epsilon _{i1}}x_{2}^{\epsilon 
_{i2}}\cdots x_{i}^{\epsilon_{ii}}, $$ 
 where $\alpha _{i},\epsilon _{ij}\in \{0,1\}$. 
 \end{defn} 

 \begin{cor} 
 Every triangular system is a fixed point system.  
 \end{cor} 

 \begin{proof} 
 Let $f$ be a triangular system with dependency graph $\X$. 
Then $\X$ consists of the glueing of components with just 
one element. Therefore each component corresponds to a 
fixed point system, and by Theorem \ref{main2} we conclude 
that $f$ is a fixed point system as well. 
 \end{proof} 

 \begin{cor} \label{componentsfps} 
 Let $f$ be a system with dependency graph $\X$.  Then $f$ 
is a fixed point system if and only if every strongly 
connected component of $\X$ either corresponds to a fixed 
point system or is connected by a walk to a zero in $\X$. 
 \end{cor} 

 \begin{proof} 
 This is an immediate consequence of Theorem \ref{main2}. 
 \end{proof} 

 {\em Remark. } 
 The order in which we enumerate the variables does not 
really matter. We will certainly get the same state space 
up to isomorphism. Namely, if 
 $$ f = (f_1(x_1, \ldots, x_n), \ldots, f_n(x_1, \ldots, 
x_n)) $$ 
 is a parallel update system and $\sigma \in S_n$ is a 
permutation, then 
 $$ \sigma f = (f_{\sigma^{-1}(1)}(x_{\sigma(1)}, \ldots, 
x_{\sigma(n)}), \ldots, f_{\sigma^{-1}(n)}(x_{\sigma(1)}, 
\ldots, x_{\sigma(n)})) $$ 
 has a state space isomorphic to the state space of $f$. In 
particular, this defines a group action of $S_n$ on the 
fixed point systems on $n$ variables. 

 \begin{thm} \label{main5} 
 Let $\X$ be the digraph of a fixed point system $f = (f_1, 
\ldots, f_n)$. If there is no walk from $a_i$ to $a_j$ or 
if $a_i$ or $a_j$ have a walk of length greater than or 
equal to $1$ to themselves, then $g := (f_1, \ldots, x_i 
f_j, \ldots, f_n)$ is a fixed point system. 
 \end{thm} 

 \begin{proof} 
 By Lemma \ref{poset} and the definition of glueing, $\X$ 
is an iterated glueing of connected components: $\X = 
(\ldots(\X_1 \# \X_2) \# \ldots ) \# \X_r$. By Corollary 
\ref{componentsfps} the connected components $\X_i$ are 
fixed point systems or they are connected by a directed 
path to a zero. The connected components that are fixed 
point systems then have loop number 0 or 1. 

 The multiplication of $f_j$ with $x_i$ introduces an extra 
edge $a_j \to a_i$ into the digraph $\X$, unless there is 
already a factor $x_i$ in $f_j$, in which case the graph 
does not change. Let $\Y$ be the digraph of $g$. We 
distinguish a number of cases. 

 {\em Case 1:} Suppose that $a_i$ and $a_j$ lie in the same 
component $\X_s$. If $i \not= j$, then $\X_s$ has loop 
number 1 or is connected to a zero. So all components 
retain their properties responsible for $\Y$ being a fixed 
point system. If $i = j$, then the loop number of $\X_s$ 
becomes 1, so $\Y$ is again a fixed point system. 

 {\em Case 2:} Let $a_i$ lie in $\X_s$ and $a_j \in \X_t$, 
and assume that there is no walk from $\X_s$ to $\X_t$. If 
there is a walk from $\X_t$ to $\X_s$ then $\X_s \geq \X_t$ 
in the partial order of $E(\X)$). The existence of an edge 
$a_j \to a_i$ does not change this property, so $\Y$ is a 
fixed point system. If there is no walk from $\X_t$ to 
$\X_s$, then we can extend the partial order of $E(\X)$ to 
a total order such that $\X_s \geq \X_t$. Then the edge 
$a_j \to a_i$ does not change the order of the glueing, so 
$\Y$ is a fixed point system. 

 {\em Case 3:} Let $a_i\in \X_s$ and $a_j \in \X_t$. Assume 
that there is a walk from $\X_s$ to $\X_t$ and thus a walk 
from $a_i$ to $a_j$. By hypothesis we have two cases. 
Either there is a closed walk $a_i \to a_i$ of length 
greater than or equal to $1$. Then the component $\X_s$ has 
loop number 1. After inserting the edge $a_j \to a_i$ the 
two components $\X_s$ and $\X_t$ are joined into one 
connected component of loop number 1. All other components 
remain unchanged. So $g$ is a fixed point system. Or there 
is a closed walk $a_j \to a_j$ of length greater than or 
equal to $1$. Then the component $\X_t$ has loop number 1. 
After inserting the edge $a_j \to a_i$ the two components 
$\X_s$ and $\X_t$ are joined again into one connected 
component of loop number 1. All other components remain 
unchanged. So $g$ is again a fixed point system. 
 \end{proof} 

 \begin{cor} 
 Let $f = (f_1, \ldots, f_n)$ be a fixed point system and 
$m$ a monomial. Then $mf = (m f_1, \ldots, m f_n)$ is a 
fixed point system. 
 \end{cor} 

 \begin{proof} It is sufficient to prove the corollary for 
the case where $m=x_i$ is a single variable. Consider first 
the system $(f_1, \ldots, x_i f_i, \ldots , f_n)$. As in 
the proof of the previous theorem we get an edge from $a_i$ 
to $a_i$ which can only change the loop number of the 
component of $a_i$ to 1. So this system is again a fixed 
point system. Now we can apply Theorem \ref{main5} $n-1$ 
times and get that $x_i f = (x_i f_1, \ldots, x_i f_n)$ is 
a fixed point system. 
 \end{proof}

 \end{document}